\documentclass[11pt]{article}

\usepackage{times}
\usepackage[utf8]{inputenc}
\usepackage[T1]{fontenc}
\usepackage{url}
\usepackage{booktabs}
\usepackage{multirow}
\usepackage{amsfonts}
\usepackage{nicefrac}
\usepackage{microtype}
\usepackage{faktor}
\usepackage{hhline}

\usepackage{thm-restate}

\usepackage{tikz}
\usepackage{mathtools, stmaryrd}
\usepackage{tocloft}
\usepackage[hidelinks]{hyperref}
\usepackage{xcolor}
\hypersetup{
    colorlinks,
    linkcolor={red!50!black},
    citecolor={blue!50!black},
    urlcolor={blue!80!black}
}
\definecolor{ao(english)}{rgb}{0.0, 0.5, 0.0}

\usepackage[algo2e, algoruled, boxed, vlined]{algorithm2e}
\usepackage{algorithm}
\usepackage{algpseudocode}

\SetKwInput{KwInput}{Input}
\SetKwInput{KwOutput}{Output}

\usepackage[noabbrev]{cleveref}

\newtheorem{Th}{Theorem}[section]

\newtheorem{Prop}[Th]{Proposition}

\newtheorem{Rem}[Th]{Remark}
\newtheorem{?}[Th]{Problem}
\newtheorem{Ex}[Th]{Example}

\newenvironment{proof}[1]
{\noindent \textit{Proof.} 
}
{
$\hfill\blacksquare$}

\def\QG{\operatorname{QG}}

\usepackage{amssymb}
\usepackage{pifont}

\usepackage{float}
\usepackage{wrapfig}

\usepackage{multicol}
\usepackage{xcolor}
\usepackage{tcolorbox}
\usepackage{amsmath}
\usepackage{mdframed}
\newmdenv[topline=false,rightline=false]{leftbot}

\DeclareMathOperator*{\argmin}{arg\,min}
\DeclareMathOperator*{\spanop}{{span}}

\newcommand{\range}[2]{\llbracket #1, #2 \rrbracket}

\newcommand{\F}{\class{\mathcal{F}}}
\newcommand{\Fml}{\class{\mathcal{F}_{\mu, L}}}
\newcommand{\Qml}{\class{\mathcal{Q}_{\mu, L}}}
\newcommand{\A}{\algo{\mathcal{A}}}

\newcommand{\xs}{x_{\star}}

\renewcommand{\leq}{\leqslant}
\renewcommand{\geq}{\geqslant}

\newcommand{\class}[1]{{\color{cyan}{#1}}}
\newcommand{\algo}[1]{{\color{purple}{#1}}}

\usepackage[letterpaper,top=2cm,bottom=2cm,left=3cm,right=3cm,marginparwidth=1.75cm]{geometry}

\title{On Fundamental Proof Structures in First-Order Optimization}

\author{
Baptiste Goujaud \thanks{Baptiste Goujaud, CMAP, Ecole Polytechnique, Institut Polytechnique de Paris, Route de Saclay, 91120 Palaiseau, France. {\tt\small baptiste.goujaud@gmail.com}.} ~ 
Aymeric Dieuleveut \thanks{Aymeric Dieuleveut, CMAP, Ecole Polytechnique, Institut Polytechnique de Paris, Route de Saclay, 91120 Palaiseau, France. {\tt\small aymeric.dieuleveut@polytechnique.edu}.} ~
Adrien Taylor \thanks{Adrien Taylor, INRIA Paris, 2 Rue Simone IFF, 75012 Paris, France. {\tt\small adrien.taylor@inria.fr}.}
}

\begin{document}

\maketitle

\begin{abstract}
    First-order optimization methods have attracted a lot of attention due to their practical success in many applications, including in machine learning. Obtaining convergence guarantees and worst-case performance certificates for first-order methods have become crucial for understanding ingredients underlying efficient methods and for developing new ones. However, obtaining, verifying, and proving such guarantees is often a tedious task. Therefore, a few approaches were proposed for rendering this task more systematic, and even partially automated. In addition to helping researchers finding convergence proofs, these tools provide insights on the general structures of such proofs. We aim at presenting those structures, showing how to build convergence guarantees for first-order optimization methods.
\end{abstract}

\section{Introduction}\label{sec:introduction}
In recent years, there has been a significant surge in the interest surrounding first-order optimization methods, primarily driven by their remarkable efficiency on a number of applications, notably within the field of machine learning~(see e.g.,~\cite{bottou2007tradeoffs}).
Theoretical foundations for those methods played a crucial role in this success, e.g., by enabling the development of momentum-type methods ~(see e.g.,~\cite{polyak_gradient_1963, nesterov1983method}). Formally, we consider the optimization problem
\begin{equation}
    x_\star \triangleq \arg\min_{x\in \mathbb R^d} f(x) \tag{OPT}\label{eq:opt}
\end{equation}
where $f$ belongs to a set $\class{\F}$ (often referred to as a ``class of functions'', e.g., the set of convex functions, the set of strongly-convex and smooth functions, or the set of quadratic convex functions, etc.). Classical first-order optimization methods for solving this problem include \emph{gradient descent}~(GD)~\cite{cauchy1847methode}, \emph{Nesterov accelerated gradient method}~(NAG)~\cite{nesterov1983method}, and the \emph{heavy-ball method}~(HB)~\cite{polyak_gradient_1963}.

In this context, a key question is to obtain a priori performance guarantees for an iterative algorithm~\algo{$\A$} (i.e., \algo{$\A$} is a rule for generating sequences of approximations $(x_t)_{t\leq T}$ to the minimizers of a certain function $f$) when the function~$f$ to be minimized belongs to a set \class{$\F$}. The most popular framework for such analyses of optimization algorithms is that of \emph{worst-case analyses}, see, e.g.,~\cite{nesterov1983method,dvurechensky2021first,bubeck2015convex,d2021acceleration,chambolle2016introduction}.
Given an algorithm \algo{$\A$}, the \emph{worst-case analysis} framework consists in finding guarantees that hold for \emph{every} function of the class.

In other words, we aim at evaluating the worst-case accuracy of \algo{$\A$} over the functions of the class~\class{$\F$} after a given number of iterations~$T$.
For doing so, there are many different possible notions of \emph{accuracy} (or performance) which we denote by $P(f, (x_t)_{t \leq T})$ and that we aim at minimizing.
Letting $x_T$ be the output of an algorithm, common examples of such metrics include the distance of the last iterate to an optimum $\|x_T-\xs\|$, the function value accuracy of the last iterate $f(x_T)-f(\xs)$, or its gradient norm $\|\nabla f(x_T)\|$.
Usually, $x_T$ can be arbitrarily bad just by choosing $x_0$ arbitrarily far away from the optimizer $x_\star$. Therefore, we usually need to assume $x_0$ to be not too bad, such as $x_0 \in \mathcal{N}(x_\star)$ where $\mathcal{N}(x_\star)$ can be any fixed set (that we call a ``neighborhood'' of the optimizer $x_\star$) and depends on $x_\star$.
Common examples of such neighborhood are balls around the optimizer $\left\{ x | \|x-x_\star\| \leq R\right\}$ or the set $\left\{ x | f(x) - f_\star \leq R\right\}$.

The smallest upper bound on $P(f, (x_t)_{t \leq T})$ that holds for any dimension $d \geq 1$, for any function~$f\in\F$, for any starting point $x_0\in\mathcal{N}(x_\star) \subset \mathbb{R}^d$, and for any $(x_t)_{t\leq T}$ generated by \algo{$\A$} applied on $f$ from $x_0$, is the optimal value to the problem of computing the worst-case:

\begin{center}
    \fbox{\parbox{.65\textwidth}{
        \begin{center}
            \begin{equation}
                \left|\hspace{-.1cm}
                \begin{array}{cc}
                    \underset{\substack{f \in\F, d \geq 1\\(x_t)_{t\leq T} \in\left(\mathbb{R}^d\right)^{T+1}}}{\text{maximize }} & P(f, (x_t)_{t\leq T}) \\
                    \text{subject to } &
                    \left\{
                    \begin{array}{c}
                        x_0 \in \mathcal{N}(x_\star) \\
                        (x_t)_{t\leq T} = \A(f, T, x_0)
                    \end{array}
                    \right.
                \end{array}
                \right.
                \tag{$\mathcal{P}$} \label{eq:problem}
            \end{equation}
        \end{center}
    }}
\end{center}

In the black-box model, iterative algorithms gather information about $f$ through so-called \emph{oracles}, which we denote by $\mathcal{O}^{(f)}$.
Classical oracles used in first-order optimization are gradient evaluations $\mathcal{O}^{(f)}(x)=\nabla f(x)$ and approximate gradients $\mathcal{O}^{(f)}(x)\approx \nabla f(x)$ (e.g., stochastic gradients), but also proximal operators (see, e.g.~\cite{combettes2011proximal}), etc.
At step $t \in \llbracket 1, T \rrbracket$, \algo{$\A$} collects oracles on the previous iterates $(\mathcal{O}^{(f)}(x_s))_{s\leq t-1}$ and outputs $x_t$ based on those information through the update function $A_t$ as $x_t = \algo{A_t}((x_s, \mathcal{O}^{(f)}(x_s))_{s<t})$.

\textbf{Notation.} For readability purposes, all notation used throughout this paper are summarized as follows.
\begin{table}[h!]
    \centering
    \begin{tabular}{ll}
\toprule
Notation & Corresponding object \\
\midrule
    $\F$ & Class of functions (generic form) \\
    $f$ & Objective function \\
    $x_\star$ & Optimal point \\

    $x_0$ & Initial iterate \\
    $\mathcal O^{(f)} $ & Generic oracle applied on $f$ \\
    $\A$ & Algorithm (generic form) \\
    $(x_t)_{t\leq T}$ & Sequence of iterates generated by $\A$, i.e. $(x_t)_{t\leq T} = \A(f, T, x_0)$ \\
    $\algo{(A_t)_{1 \leq t \leq T}}$ & Update function of the algorithm $\A$, i.e. $\forall t, x_t = \algo{A_t}((x_s, \mathcal{O}^{f}(x_s))_{s<t})$ \\
    $T$ & Total number of iterations \\
    $t$ & Current iteration index \\

    $\Fml$ & Class of $L$-smooth and $\mu$-strongly-convex functions ($0\leq \mu\leq L$)\\
    $\Qml$ & Class of $L$-smooth and $\mu$-strongly  convex quadratic functions ($0\leq \mu\leq L$)\\

    $(V_t)_t$ & Lyapunov sequence \\
    $F, G$  & Linearization variables (after SDP lifting)\\
    $P(f, (x_t)_{t\leq T})$ & Performance metric \\
    \bottomrule\\[-.5cm]
\end{tabular}
\end{table}

\textbf{Outline.} In \Cref{sec:explicit_implicit_class}, we discuss two ways of characterizing classes of functions and detail the main cases for which we can solve~\eqref{eq:problem}. In \Cref{sec:explicit_implicit_algo}, we discuss an alternative way of describing the algorithm \algo{$\A$} simplifying the resolution of~\eqref{eq:problem}. \Cref{sec:proof_structures} outlines a systematic approach for acquiring proofs of worst-case performance certificates and delves into their underlying structures. We further elaborate on how this structure can be exploited for extending the applicability range of the worst-case guarantees. Among others, we show how the properties of these proofs allow building algorithms. Finally, \Cref{sec:lyap} provides a natural approach for discovering Lyapunov sequences.

\section{From explicit to implicit classes of functions}\label{sec:explicit_implicit_class}

This section describes two ways of specifying a class of functions as part of the worst-case analysis of a given algorithm.
We describe two different methods to approach and solve~\eqref{eq:problem} depending on the ways $\F$ is specified.
More specifically, we focus on two specific classes of functions to illustrate our explanations, namely $L$-smooth $\mu$-strongly convex quadratic functions (notation $\Qml$) and $L$-smooth $\mu$-strongly convex functions (notation $\Fml$).

\subsection{Convex quadratic optimization}\label{subsec:quad}
First-order optimization methods were extensively studied in the context of minimizing quadratic convex functions.
Such functions can be described \textbf{explicitly} as

\begin{equation}
    f(x) \triangleq \frac{1}{2}(x - x_\star)^T H (x - x_\star) + f_\star,
\end{equation}
where $H$ is the symmetric positive semi-definite Hessian of~$f$, $x_\star$ its optimizer and $f_\star$ its minimal value.
This expression allows to explicitly compute the gradient $\nabla f(x) = H(x - x_\star)$, and first-order optimization methods can be expressed through polynomials due to the following property~(e.g.,~\cite[Prop.4.1]{goujaud2022super}).

\begin{leftbot}
    \vspace{-.2cm}
    \begin{Prop}
    \label{prop:link_algo_poly}
        Let $f\in\class{\mathcal{Q}_{0, \infty}}$ and $x_0\in\mathbb{R}^d$.
        It holds that
        \vspace{-.2cm}
        \begin{equation}
            x_{t+1} \in x_0 + \spanop\{ \nabla f(x_0),\ldots, \nabla f(x_t) \} \,, \label{def:first_order_algo}
        \end{equation}
        if and only if there exists a sequence of polynomials $(P_t)_{t\in\mathbb{N}}$, each of degree at most 1 more than the highest degree of all previous polynomials and $P_0$ of degree 0 (hence the degree of $P_t$ is at most $t$), such that
        \begin{equation}
            \vphantom{\sum_i}\forall \;t \quad x_t - x_\star = P_t(H)(x_0-x_\star), \quad P_t(0)=1\,. \label{eq:link_polynomial}
        \end{equation}
    \end{Prop}
\end{leftbot}

In this context,~\eqref{eq:problem} can be solved by solving a polynomial problem of the form $\max_{H} \|P_t(H)\|$ where $H$ is a symmetric matrix verifying some conditions (e.g. $\mu I \preceq H \preceq L I$ when $f\in\class{\mathcal{Q}_{\mu, L}}$).
This link between first-order algorithms and polynomials has been used by~\cite{golub1961chebyshev} for discovering the Chebyshev method and by~\cite{polyak_gradient_1963} for the ``heavy-ball'' method, still used nowadays far beyond quadratic optimization (e.g.\ in stochastic optimization of neural networks~\cite{sutskever2013importance}).
This property has also been exploited more recently for obtaining new algorithms with provable guarantees on quadratic functions~(see e.g.,~\cite{fischer2011polynomial, scieur2018acceleration, d2021acceleration, pedregosa2020acceleration, scieur2020universal, berthier2020accelerated, goujaud2022super, goujaud2022quadratic, cunha2022only}).

\subsection{Infinite-dimensional spaces of functions}\label{subsec:implicit_class}

As opposed to previous sections, many classes of functions are described \emph{implicitly} as regions  of infinite-dimensional spaces of functions.
In other words, such functions are defined by sets of inequalities.
This section deals with the analyses of such classes.
This is due to the fact the set of all functions of the class are not described by a finite number of parameters, but rather by constraints (inequalities).
Studying~\eqref{eq:problem} for classes that are defined \textbf{implicitly} through sets of constraints appears to be much less natural.
In this situation,~\eqref{eq:problem} is often referred to as a \emph{performance estimation problem (PEP)}~\cite{drori2014performance, taylor2017smooth, taylor2017exact}.
This tool primarily relies on two crucial components: interpolation conditions and SDP lifting.

\textbf{Interpolation conditions.}
We remark that the description of the algorithm and the objective of~\eqref{eq:problem} both only depend on
the oracle values of $f$ on the iterates $(x_{t})_{t \leq T}$.
We introduce the variables $(\mathcal{O}_t)_{t\leq T}$.
The constraint $f\in\F$ must be replaced by the constraint that there exists at least one element $f \in \F$ such that $(\mathcal{O}^{(f)}(x_t))_{t\leq T} = (\mathcal{O}_t)_{t\leq T}$ ($\mathcal{O}_t$ is a reachable value for $\mathcal{O}^{(f)}(x_t)$, when $f\in\F$). As an example, $f_t$ and $g_t$ are potential values of respectively $f(x_t)$ and $\nabla f(x_t)$.
Formally, we define the equivalence relation $\sim_{\eqref{eq:problem}}$ as $f_1 \sim_{\eqref{eq:problem}} f_2$ if and only if $\forall t \in\range{0}{T}\cup\{\star\}, \mathcal{O}^{(f_1)}(x_t) = \mathcal{O}^{(f_2)}(x_t)$.
Since the only information $\A$ gathers on $f$ is the oracle outputs at the iterates $x_t$,
two functions coming from the same equivalence class both produce feasible points of~\eqref{eq:problem} with the same objective value.
In other words, those two functions are undistinguishable using only the information available to $\A$.
We can therefore rewrite~\eqref{eq:problem} in terms of $(\mathcal{O}_t)_{t\leq T} \in \class{\F / \sim_{\eqref{eq:problem}}}$ instead of $f \in \F$, so that the set of optimization variables now lives in finite dimension.
This constraint is referred to as \emph{interpolation conditions}.

\begin{leftbot}
\begin{Ex}[First-order algorithm on $\Fml$]
    \label{ex:fom_on_strooth}
    Let $L\geq \mu>0$ two positive real numbers.
    A function $f$ is $L$-smooth and $\mu$-strongly-convex when $f$ is continuously differentiable and verifies the two inequalities:
    \begin{align}
        f(x) & \leq f(y) + \left<\nabla f(y), x - y\right> + \frac{L}{2} \| x - y \|^2, \\
        f(x) & \geq f(y) + \left<\nabla f(y), x - y\right> + \frac{\mu}{2} \| x - y \|^2,
    \end{align}
    for all $x, y$ and where $\nabla f$ denotes the gradient of $f$.

    Studying a first-order algorithm (i.e.\ an algorithm based on the oracle $\mathcal{O}^{(f)} \triangleq (\nabla f, f)$) on the class $\Fml$ appears to be challenging at first sight due to the infinite number of parameters needed for describing $\Fml$.
    However,~\cite[Theorem 4]{taylor2017smooth} provides interpolation conditions for the class $\Fml$ of $L$-smooth $\mu$-strongly-convex functions and enables an exact study of the worst-case of several algorithms on this class of functions:
\begin{align}
    \forall i, j,~ f_i \geq & f_j + \left< g_j, x_i - x_j \right> + \tfrac{1}{2L}\|g_i - g_j\|^2 \tag{IC} \label{eq:cni} \\
    & + \tfrac{\mu}{2(1 - \mu / L)}\|x_i - \tfrac{1}{L}g_i - x_j + \tfrac{1}{L}g_j\|^2. \nonumber
\end{align}

Indeed, in this case,~\eqref{eq:problem} can be written in finite dimension as
\begin{equation*}
    \left|
    \begin{array}{l}
        \underset{(g_t, f_t)_{t\leq T}\in\left(\mathbb{R}^d \times \mathbb{R}\right)^{T+1}}{\underset{d \geq 1, (x_t)_{t\leq T} \in\left(\mathbb{R}^d\right)^{T+1}, x_\star\in\mathbb{R}^d, }{\text{maximize }}} P((x_t, g_t, f_t)_{t\leq T}) \\
        \text{s.t.}
        \left\{
        \begin{array}{rl}
            & \hspace{-.5cm} x_0 \in \mathcal{N}(x_\star) \\
            & \hspace{-.5cm} \forall t \leq T, ~ x_t = \algo{A_t}((x_s, \mathcal{O}^{f}(x_s))_{s<t}) \\
            & \hspace{-.5cm} \forall i, j, ~ \class{f_i \geq f_j + \left< g_j, x_i - x_j \right> + \tfrac{1}{2L}\|g_i - g_j\|^2} \\
            & \hspace{1.2cm} \class{+ \tfrac{\mu}{2(1 - \mu / L)}\|x_i - \tfrac{1}{L}g_i - x_j + \tfrac{1}{L}g_j\|^2.}
        \end{array}
        \right.
    \end{array}
    \right.
\end{equation*}

\end{Ex}
\end{leftbot}

\textbf{SDP lifting.}
In many cases (see, e.g, \Cref{ex:fom_on_strooth}, and \cite[Theorem 3.5]{taylor2017exact}), interpolation conditions are written in terms of quadratic and bilinear expressions of $x_t$ and $g_t$ and linear expressions of $f_t$.
Because of the quadratic dependency in $x_t$ and $g_t$, this problem is generally non-convex.
\emph{SDP lifting} can convexify this problem if all other parts of this problem also contain only quadratic expressions of $x_t$ and $g_t$.
For example, classical choices for $P(f, (x_t)_{t\leq T})$ are $\|x_T-\xs\|^2$, $f(x_T)-f(\xs)$, or $\|\nabla f(x_T)\|^2$.
Similarly, a classical choice for $x_0 \in \mathcal{N}(x_\star)$ is $\|x_0 - x_\star\|^2 \leq R^2$ for some radius $R>0$.
Finally, the updates $(A_t)_t$ of the algorithm $\A$ are often of the form

\begin{equation}
    x_{t} = A_t((x_s, \nabla f(x_s), f(x_s))_{s\leq t-1}) = x_0 - \sum_{s=0}^{t-1} \gamma_s^{(t)} \nabla f(x_s) \label{eq:fom}
\end{equation}
for some sequence of scalars $(\gamma_s^{(t)})_{s\in \llbracket 0, t-1 \rrbracket}$.
Substituting $x_t$ for $t\geq 1$ in the problem by their corresponding expressions given by~\eqref{eq:fom} preserves the above observation: the dependency of~\eqref{eq:problem} in $(x_t, g_t)_{t\leq T}$ is exclusively quadratic.
Actually, in this specific case, all occurrences of $(x_t)_{t\geq 1}$ have been replaced by linear combinations of $x_0$ and $(g_t)_{t\leq T}$.
\emph{SDP lifting} consists in introducing the Gram matrix $G$ of $(x_0 - x_\star, (g_t)_{t\leq T})$.
This way, all quadratic expressions of $(x_t, g_t)_{t\leq T}$ are linear combinations of the entries of $G$.
We also introduce the vector $F$ storing the values $(f_t - f_\star)_{t\leq T}$.

Finally~\eqref{eq:problem} is rewritten with linear objective and constraints only as well as an SDP constraint $G \succeq 0$.
\begin{equation*}
    \left|
    \begin{array}{cc}
        \underset{F, G \succeq 0}{\text{maximize }} & ~ \left<F, v_{P}\right> + \left<G, M_{P}\right> \\
        \text{subject to } &
        \left\{
        \begin{array}{rc}
            & \left<F, v_{I}\right> + \left<G, M_{I}\right> \leq R^2 \\
            \forall k, & \hspace{-.2cm} \left<F, \class{v_{\F}^{(k)}}\right> + \left<G, \class{M_{\F}^{(k)}}\right> \leq 0
        \end{array}
        \right.
    \end{array}
    \right.
\end{equation*}
Vectors $(v_P, v_I, \class{(v_{\F}^{(k)})_k})$ and matrices $(M_P, M_I, \class{(M_{\F}^{(k)})_k})$ are constants depending on the algorithm $\A$, the class $\F$, and the performance metric $P$ under consideration.
More specifically, indices $P$, $I$ and $\F$ respectively correspond to the performance metric, the initialization constraint and the class interpolation conditions.
The algorithm is directly encoded in the fact that $G$ does not contain inner product with $(x_t)_{t\geq 1}$.
As an example, to express $\|x_T - x_\star\|^2$ in terms of $G$, one needs to actually choose $M$ with $\left<G, M\right> = \|x_0 - x_\star - \sum_{s=0}^{T-1} \gamma_s^{(T)} \nabla f(x_s)\|^2$.

\textbf{Key conditions.} The above procedure generally works under the following conditions:
\begin{itemize}
    \item $\A$ is a first-order algorithm whose updates $\algo{(A_t)_t}$ can be expressed linearly in terms of observed gradients;
    \item The interpolation constraints of the class of functions $\F$ are known and expressible linearly in $F$ and $G$;
    \item The performance metric as well as the initial condition are also expressible linearly in terms of $F$ and $G$.
\end{itemize}

Many pairs of function class and algorithm meet the right conditions and have been studied using the PEP framework.
Tools in Matlab~\cite{taylor2017performance} and Python~\cite{goujaud2022pepit} have been implemented to automate this task and provide worst-case guarantees.
Many examples of usages are listed in the corresponding documentations.

\section{From explicit to implicit algorithms}\label{sec:explicit_implicit_algo}

So far, we only considered explicit algorithms of the form~\eqref{eq:fom}.
Note that, just as for classes of functions, algorithms can be expressed implicitly via sets of (in)equalities.
This is the case for line-search based algorithms.
Indeed, the step-size associated with line-search is not uniform over the problem class, therefore algorithms containing line-search update cannot be written as~\eqref{eq:fom}, and therefore do not meet the key conditions mentioned in the previous section.
A relaxation of the gradient descent with exact line-search has been proposed in~\cite{de2017worst}.
Since this algorithm cannot be written as~\eqref{eq:fom} with pre-determined $\gamma_s^{(t)}$, we cannot specify $(x_t)_{t\geq 1}$ in terms of $x_0$ and $(g_t)_{t\leq T}$.
Therefore, all vectors $(x_t, g_t)_{t\leq T}$ must be considered as linearly independent.
For this problem, $G$ is the Gram matrix of all $(x_t, g_t)_{t\leq T}$.

Therefore, the algorithm is not totally encoded in $v_P, M_P, v_I, M_I, \class{v_{\F}}$ and $\class{M_{\F}}$ anymore and must be specified by new constraints.
In particular, the updates of gradient descent with line-search verify that
\begin{align}
    \left< g_{t+1}, g_t \right> = 0 \\
    \left< g_{t+1}, x_{t+1} - x_t \right> = 0
\end{align}
As for all the other elements of~\eqref{eq:problem}, those constraints only involve quadratic terms of $(x_t, g_t)_{t\leq T}$ and can therefore be expressed linearly in terms of~$G$, parametrized by the vectors~$\algo{(v_{\A}^{(l)})_l}$ and matrices~$\algo{(M_{\A}^{(l)})_l}$. This time,~\eqref{eq:problem} writes
\begin{equation}
    \left|
    \begin{array}{cc}
        \underset{F, G \succeq 0}{\text{maximize }} & \hspace{.2cm} \left<F, v_{P}\right> + \left<G, M_{P}\right> \\
        \text{subject to } &
        \left\{
        \begin{array}{rcl}
            & \left<F, v_{I}\right> + \left<G, M_{I}\right> & \leq R^2 \\
            \forall k, & \hspace{-.2cm} \left<F, \class{v_{\F}^{(k)}}\right> + \left<G, \class{M_{\F}^{(k)}}\right> & \leq 0 \\
            \forall l, & \hspace{-.2cm} \left<F, \algo{v_{\A}^{(l)}}\right> + \left<G, \algo{M_{\A}^{(l)}}\right> & \leq 0 \\
        \end{array}
        \right.
    \end{array}
    \right.
    \tag{PEP-primal} \label{eq:pep_primal}
\end{equation}

\section{Proof structures in first-order optimization}\label{sec:proof_structures}

There is an extensive literature on first-order optimization, offering a broad range of possibly advanced worst-case guarantees and their associated proofs.
In the previous sections, we saw conditions under which the problem of computing worst-case guarantees was tractable.
In this section, we detail how to obtain proofs from PEPs and what we can conclude on the general structure of proofs in first-order optimization.

\subsection{Obtaining proofs with PEPs}\label{subsec:proofs_with_pep}

Thanks to interpolation conditions and SDP lifting,~\eqref{eq:problem} rewrites as a convex optimization problem. We consider the dual of the problem. 
Let's then introduce the Lagrangian multipliers $\tau$, $\class{(\lambda_{\F}^{(k)})_k}$, $\algo{(\lambda_{\A}^{(l)})_l}$ associated to the constraints of~\eqref{eq:pep_primal}.

\vspace{-.1cm}
\begin{equation*}
    \left|
    \begin{array}{cc}
        \underset{F, G \succeq 0}{\text{maximize }} & \hspace{-1.1cm} \left<F, v_{P}\right> + \left<G, M_{P}\right> \\
        \text{subject to } &
        \left\{
        \begin{array}{rcll}
            & \hspace{-.5cm} \left<F, v_{I}\right> + \left<G, M_{I}\right> & \hspace{-.2cm} \leq R^2 & \hspace{-.2cm} : \tau \\
            \forall k,~ & \hspace{-.5cm} \left<F, \class{v_{\F}^{(k)}}\right> + \left<G, \class{M_{\F}^{(k)}}\right> & \hspace{-.2cm} \leq 0 & \hspace{-.2cm}: \class{\lambda_{\F}^{(k)}} \\
            \forall l,~ & \hspace{-.5cm} \left<F, \algo{v_{\A}^{(l)}}\right> + \left<G, \algo{M_{\A}^{(l)}}\right> & \hspace{-.2cm} \leq 0 & \hspace{-.2cm}: \algo{\lambda_{\A}^{(l)}} \\
        \end{array}
        \right.
    \end{array}
    \right.
\end{equation*}
\noindent
The Lagrangian then writes
\begin{align*}
    \mathcal{L} \triangleq &
    \left<F, v_{P}\right> + \left<G, M_{P}\right> - \tau \left[ \left<F, v_{I}\right> + \left<G, M_{I}\right> - R^2 \right] \\
    & - \sum_{k} \class{\lambda_{\F}^{(k)}} \left[ \left<F, \class{v_{\F}^{(k)}}\right> + \left<G, \class{M_{\F}^{(k)}}\right> \right] \\ & - \sum_l \algo{\lambda_{\A}^{(l)}} \left[ \left<F, \algo{v_{\A}^{(l)}}\right> + \left<G, \algo{M_{\A}^{(l)}}\right> \right] \\
    = &
    \tau R^2 + \left<F, v_{P} - \tau v_{I} - \sum_{k} \class{\lambda_{\F}^{(k)}} \class{v_{\F}^{(k)}} - \sum_l \algo{\lambda_{\A}^{(l)}} \algo{v_{\A}^{(l)}} \right> \\
    & + \left<G, M_{P} - \tau M_{I} - \sum_{k} \class{\lambda_{\F}^{(k)}} \class{M_{\F}^{(k)}} - \sum_l \algo{\lambda_{\A}^{(l)}} \algo{M_{\A}^{(l)}}\right>
\end{align*}
\noindent
The dual is obtained by maximizing over the primal variables:
\vspace{-.2cm}
\begin{equation}
    \left|
    \begin{array}{l}
        \underset{\tau, \class{\lambda_{\F}^{(k)}}, \algo{\lambda_{\A}^{(l)}} \geq 0}{\text{minimize }} \tau R^2  \\
        \text{s.t.}
        \left\{
        \begin{array}{lc}
            v_{P}\  - \tau v_{I}\ \,  - \sum_{k} \class{\lambda_{\F}^{(k)}} \class{v_{\F}^{(k)}} \ \, - \sum_l \algo{\lambda_{\A}^{(l)}} \algo{v_{\A}^{(l)}} &  = 0 \\
            M_{P} - \tau M_{I} - \sum_{k} \class{\lambda_{\F}^{(k)} M_{\F}^{(k)}} - \sum_l \algo{\lambda_{\A}^{(l)} M_{\A}^{(l)}} &  \preceq 0
        \end{array}
        \right.
    \end{array}
    \right. \label{eq:pep_dual}
        \tag{PEP-dual} 
\end{equation}
\noindent
For any feasible primal $F, G$ and feasible dual $\tau, \class{(\lambda_{\F}^{(k)})_{k}}, \algo{(\lambda_{\A}^{(l)})_l}$, we know the objective of the dual is larger than the Lagrangian value, that is:

\begin{align}
    & \underbrace{\left<F, v_{P}\right> + \left<G, M_{P}\right>}_{\text{Performance metric}} - \tau \underbrace{\left[ \left<F, v_{I}\right> + \left<G, M_{I}\right> \right]}_{\text{Initialization}} \nonumber \\
    & \leq \sum_{k} \class{\lambda_{\F}^{(k)}} \underbrace{\left[ \left<F, \class{v_{\F}^{(k)}}\right> + \left<G, \class{M_{\F}^{(k)}}\right> \right]}_{\text{Class constraint}} \nonumber \\
    & \quad + \sum_l \algo{\lambda_{\A}^{(l)}} \underbrace{\left[ \left<F, \algo{v_{\A}^{(l)}}\right> + \left<G, \algo{M_{\A}^{(l)}}\right> \right]}_{\text{Algorithm constraint}} \nonumber \\
    & \leq 0.
    \tag{Generic proof} \label{eq:generic_proof}
\end{align}

\noindent
In words, the proof of a worst-case guarantee is obtained by  linearly combining  all  available constraints,  with coefficients that are the dual variables of the PEP.
Indeed, the difference between the performance metric and $\tau$ times the initialisation measure of proximity to the optimizer is decomposed as the sum of three terms. The two first ones respectively correspond to the values that are enforced to be negative by the class of functions and the algorithm. The third one is called the residual and is the opposite of a sum of squares of iterates and gradients. An example of full derivation of such a proof is provided in~\Cref{sec:example}.

\begin{Rem}[No duality gap]
There generally exists a feasible point $G, F$ with $G \succ 0$, i.e.\ verifying the Slater's condition~(see~\cite{slater1950lagrange}), therefore guaranteeing strong duality of the convex reformulation of~\eqref{eq:problem}.
To ensure this, one needs to carefully remove iterates $x_t$ from the basis of $G$ when $x_t$ is completely identified from other vectors.
For instance, leaving $x_1$ in the basis of $G$ with the constraint $\|x_1 - (x_0 - \gamma g_0)\|^2 = 0$ instead of replacing $x_1$ by $x_0 - \gamma g_0$ everywhere, creates an empty interior and can break strong duality.
Each time there is no feasible $G$ with $G\succ 0$, we conclude that there is a linear relationship between elements of the basis $G$ is the Gram matrix of.
Therefore, maximally reducing the dimension of $G$ ensures strong duality.
\end{Rem}

\subsection{Understanding proofs with PEPs}\label{subsec:analyzing_proofs}

Obtaining dual feasible points provides valuable insights into essential aspects pertaining to both the class of functions under consideration and the algorithm employed to achieve the associated worst-case guarantee.

\textbf{Extension to broader sets of algorithms.}
\cite{drori2020efficient} exploit these insights to \emph{design worst-case optimal algorithms}. The authors' key observation is that~\eqref{eq:generic_proof} does not rely on all constraints to hold, but rather only on a linear combination of them. Therefore, if instead of assuming that, $\forall l, \left<F, \algo{v_{\A}^{(l)}}\right> + \left<G, \algo{M_{\A}^{(l)}}\right> \leq 0$, we can simply assume that $\sum_l \algo{\lambda_{\A}^{(l)}} \left[ \left<F, \algo{v_{\A}^{(l)}}\right> + \left<G, \algo{M_{\A}^{(l)}}\right> \right] \leq 0$, therefore relaxing a lot of assumptions about the algorithm and then generalizing the proof to all the algorithms verifying the remaining assumption. This was applied to the impractical algorithm~\eqref{eq:gfom} described as follow:
\begin{equation}
    \forall t, x_{t+1} = \underset{x\in x_0 + \spanop \left\{ \nabla f(x_0), \cdots, \nabla f(x_t) \right\}}{\argmin} f(x), \tag{GFOM} \label{eq:gfom}
\end{equation}
greedily minimizing the objective value in the affine space of all the observed directions.
For some classes of functions, this algorithm is worst-case optimal. This is the case, for instance, for the class of quadratic convex functions on which~\eqref{eq:gfom} is equivalent to the so-called \emph{conjugate gradient method}. This is also the case for the class of $L$-smooth convex functions, allowing to find a broad range of worst-case optimal algorithms on this class, including the so-called \emph{optimized gradient method (OGM)} \cite{drori2017exact, drori2020efficient}. Generating such worst-case optimal algorithms works as follow:
\begin{enumerate}
    \item We note that~\eqref{eq:gfom} verifies the following orthogonality constraints:
    \begin{equation}
        \forall t, \left\{
        \begin{array}{c}
            \forall s < t, \left< g_t, g_s \right> = 0, \\
            \forall s \leq t, \left< g_t, x_s - x_0 \right> = 0.
        \end{array}
        \right. \label{eq:constraints_grouping}
    \end{equation}
    Note that following those constraints does not necessarily imply that~\eqref{eq:pep_primal}'s primal variables optimal values describe~\eqref{eq:gfom}. Nevertheless, a sufficient condition on the class $\F$ under consideration for that to happen is that
    $\F$ is contraction-preserving~(see~\cite[Definition 3]{drori2020efficient}), which happens to be the case for $\Fml$ for example.
    \item We call the corresponding dual variables $(\beta_{t, s})_{s<t}$ and $(\gamma_{t, s})_{s\leq t}$ and collect their optimal values $(\beta_{t, s}^{\star})_{s<t}$ and $(\gamma_{t, s}^{\star})_{s\leq t}$: it happens that those values can be obtained in closed-form.
    \item We group all the constraints as in~\eqref{eq:constraints_grouping}, and conclude that the worst-case guarantee of~\eqref{eq:gfom}, as well as the corresponding proof, would hold if
    \begin{equation}
        \forall t, \left< g_t, \sum_{s=0}^{t-1} \beta_{t, s} g_s + \sum_{s=0}^t \gamma_{t, s} (x_s - x_0) \right> \leq 0.
    \end{equation}
    \item When $\gamma_{t, t} \neq 0$, we conclude that, in particular, the algorithm described by the iteration
    \begin{equation}
        \forall t, x_t = x_0 - \sum_{s=0}^{t-1} \tfrac{\gamma_{t, s}}{\gamma_{t, t}} (x_s - x_0) + \frac{\beta_{t, s}}{\gamma_{t, t}} g_s
    \end{equation}
    annihilates the vector in the right-hand position of the inner product. Therefore, the worst-case guarantee of~\eqref{eq:gfom} also applies to~$\A$, using the exact same proof.
\end{enumerate}

This method has more recently been used in~\cite[Th.2.4-Cor.2.5]{goujaud2022optimal} to derive the worst-case optimal algorithm
\begin{equation}
    x_t = \frac{t}{t+1}x_{t-1} + \frac{1}{t+1}x_0 - \frac{1}{t+1}\sum_{s=0}^{t-1}\frac{1}{L} g_s \tag{HB} \label{eq:hb_qg}
\end{equation}
under the class of convex and $L$-quadratically upper bounded (\class{$L$-$\QG^+$}) functions.

\textbf{Extension to broader classes of functions.}
Interestingly, \eqref{eq:hb_qg} was studied several years ago in~\cite{ghadimi2014global} on the class $\class{\F_{0, L}}$ of $L$-smooth convex functions, itself included in the class of $L$-$\QG^+$ convex functions. On the other hand, the obtained guarantee was not better on $\class{\F_{0, L}}$ than the one obtained on the class of \class{$L$-$\QG^+$} convex functions. This shows that the guarantee obtained on $\class{\F_{0, L}}$ can be obtained using only the interpolation constraints of the class of \class{$L$-$\QG^+$} convex functions, which is a subset of the set of interpolation constraints of \class{$\F_{0, L}$}. In general, for a given class and a given algorithm, when $\class{\lambda_{\F}^{(k)}}=0$ in~\eqref{eq:generic_proof}, we conclude that the corresponding constraint has not been used. This allows to discard all the useless constraints and the result naturally holds on a larger class of functions.

\textbf{Fewer class constraints allows new algorithms.}
Most of the time, we study a family of classes of functions, parametrized by some value $L$. A classical example of this is the class of $L$-smooth convex functions $\class{\F_{0, L}}$.
The underlying \emph{interpolation constraints} $\left<F, \class{v_{\F}^{(k)}(L)}\right> + \left<G, \class{M_{\F}^{(k)}(L)}\right> \leq  0$ then depend on $L$. We generally derive and study an algorithm on $\class{\F_{0, L}}$, and obtain a guarantee that holds for any $L$ such that $\left<F, v_{P}(L)\right> + \left<G, M_{P}(L)\right> - \tau(L) \left[ \left<F, v_{I}(L)\right> + \left<G, M_{I}(L)\right> \right] \leq 0$.
The underlying algorithm can (and usually does) therefore depend on this value that is sometimes hard to access in practice. Using line-search steps is a way to get rid of the dependence on $L$ (there exists for instance line-search version of \emph{OGM} and~\eqref{eq:hb_qg} that do not involve $L$), but an exact line-search step is often not available neither. On the other hand, \emph{backtracking line-search} have been proposed~\cite{armijo1966minimization} to replace the class parameter $L$ by any surrogate value $\hat{L}$ that validates all the inequalities that are used. Indeed, we know that for any $L$,
\begin{align}
    & \underbrace{\left<F, v_{P}(L)\right> + \left<G, M_{P}(L)\right>}_{\text{Performance metric}} - \tau(L) \underbrace{\left[ \left<F, v_{I}(L)\right> + \left<G, M_{I}(L)\right> \right]}_{\text{Initialization}} \nonumber\\
    & \leq
    \sum_{j} \lambda^{(j)} \underbrace{\left[ \left<F, v^{(j)}(L)\right> + \left<G, M^{(j)}(L)\right> \right]}_{\text{Constraint}} \nonumber\\
    & \leq 0
\end{align}

Therefore, even if we do not have access to $L$, being able to find some $\hat{L}$ in an online manner such that all the surrogate constraints $\left<F, \class{v_{\F}^{(k)}(\hat{L})}\right> + \left<G, \class{M_{\F}^{(k)}(\hat{L})}\right> \leq  0$ hold, allows tuning the algorithm online with this $\hat{L}$ and obtain the guarantee 
\begin{equation*}
    \left<F, v_{P}(\hat{L})\right> + \left<G, M_{P}(\hat{L})\right> - \tau(\hat{L}) \left[ \left<F, v_{I}(\hat{L})\right> + \left<G, M_{I}(\hat{L})\right> \right] \leq 0.
\end{equation*}
We would like to apply bisection search to find such $\hat{L}$, and all we need for that is being able to verify the constraints $\left<F, \class{v_{\F}^{(k)}(\hat{L})}\right> + \left<G, \class{M_{\F}^{(k)}(\hat{L})}\right> \leq 0$ online. Note however that some constraints may involve the optimizer $x_\star$ or the minimal value $f_\star$ and are then not verifiable. The authors of~\cite[Remark 4.9]{d2021acceleration} and~\cite{park2021optimal} discuss this issue. They note that we only need to verify constraint that actually involve $L$ and that the ones that are problematic are the ones that involve both $L$ and an unknown value. They conclude that, if the dual values associated with these problematic constraints are set to 0, they are not used, and then we can proceed to \emph{backtracking line-search}. They also enforce it by removing those inequalities (or lowering them) and searching for methods that holds on this larger class of functions (verifying less inequalities) in order to be able to apply \emph{backtracking line-search} to get rid of the requirement of knowledge of the parameter class.

\section{Example: gradient descent with exact line-search}\label{sec:example}

For sake of better comprehension of the formal reasoning made in~\Cref{subsec:implicit_class,sec:explicit_implicit_algo,sec:proof_structures},
we detail in this section the development of a proof of convergence guarantee of the form~\eqref{eq:generic_proof} on an example: the gradient descent method with exact line-search, defined as
\begin{equation}
    \forall t\in\range{1}{T}, ~ x_{t} = \argmin_{x \in x_{t-1}+\spanop\left\{\nabla f(x_{t-1}) \right\}} f(x). \tag{GDLS} \label{eq:gdls}
\end{equation}
More precisely, we chose to consider the function value as performance metric, and therefore seek for a guarantee of the form
\begin{equation}
    f(x_1) - f_\star \leq \tau (f_0 - f_\star),
\end{equation}
with an appropriate $\tau$.
Note this problem has been solved in~\cite[Theorem 1.2]{de2017worst}. Here we detail how to find such a guarantee and its proof in a very systematic way, relying on the framework presented in the present tutorial.

The problem can therefore be summarized as follow:
\begin{itemize}
    \item The objective function belongs to the class $\Fml$ of $L$-smooth $\mu$-strongly-convex functions, i.e. verifies the interpolation constraints~\eqref{eq:cni},
    \item We have access to the oracle $\mathcal{O}^f(x)$ verifying:
    \begin{itemize}
        \item $\mathcal{O}^f(x)\in x + \spanop\left\{\nabla f(x)\right\}$,
        \item $\left<\nabla f\left(\mathcal{O}^f(x)\right), \nabla f(x)\right>=0$,
    \end{itemize}
    \item The algorithm $\A$ iteratively computes the update $x_{t} = \algo{A_t}((x_s, \mathcal{O}^{f}(x_s))_{s<t}) \triangleq \mathcal{O}^f(x_{t-1})$,
    \item We study exactly one step of this algorithm. That is, we want a guarantee on $x_1$ given $x_0$.
    \item The performance metric that we use is the function value $f(x_1) - f_\star$.
    \item The neighborhood $\mathcal{N}(x_\star)$ we assume $x_0$ belongs to is also define by the function value as $\left\{ x ~|~ f(x) - f_\star \leq R^2 \right\}$ for some positive $R$.
\end{itemize}

In summary, the problem~\eqref{eq:problem} writes
\begin{center}
    \fbox{\parbox{0.65\textwidth}{
        \begin{center}
            \begin{equation}
                \left|\hspace{-.1cm}
                \begin{array}{cc}
                    \underset{\substack{f \in\Fml, d \geq 1\\(x_\star, x_0, x_1)\in\left(\mathbb{R}^d\right)^{3}}}{\text{maximize }} & f(x_1) - f_\star \\
                    \text{subject to } &
                    \left\{
                    \begin{array}{c}
                        f(x_0) - f_\star \leq R^2 \\
                        (x_t)_{t \leq 1} = \text{\ref{eq:gdls}}(f, T=1, x_0)
                    \end{array}
                    \right.
                \end{array}
                \right.
            \end{equation}
        \end{center}
    }}
\end{center}

\ref{eq:gdls}'s update is defined through an optimization problem. Implementing it into the PEP framework is not straightforward. Instead, we replace the strict definition of the update by first order optimality conditions of the line search procedure:
\begin{align*}
    \algo{\left< \nabla f(x_1), \nabla f(x_0) \right>} & \algo{= 0,} \\
    \algo{\left< \nabla f(x_1), x_1 - x_0 \right>} & \algo{= 0.}
\end{align*}

Note the second one is verified because $x_1-x_0$ is colinear with $g_0$ and therefore those 2 conditions seem redundant.
However, removing the proper definition of~\eqref{eq:gdls} makes $x_1-x_0$ and $g_0$ non-necessarily colinear anymore, and the two orthogonality conditions are complementary.

Note furthermore that, replacing the actual definition of~\eqref{eq:gdls} by some conditions the latter verifies leads to a guarantee that holds over all the algorithms that verify those conditions.
This is therefore possibly a relaxation, but the result still holds.
Moreover, in this special case, and because we used the two orthogonality conditions and not just one, replacing the definition of~\eqref{eq:gdls} by those conditions is tight.
This technical assertion is based on the fact the class $\Fml$ is \emph{contraction-preserving}. This reasoning is detailed in~\cite{drori2020efficient}.

Expressing the constraints of the algorithm and the class, we obtain
\begin{equation*}
    \left|
    \begin{array}{l}
        \underset{(g_0, g_1)\in\left(\mathbb{R}^d\right)^{2}, ~ (f_\star, f_0, f_1) \in \mathbb{R}^3}{\underset{d \geq 1, (x_\star, x_0, x_1)) \in\left(\mathbb{R}^d\right)^{3}, }{\text{maximize }}} f(x_1) - f_\star \\
        \text{s.t.}
        \left\{
        \begin{array}{ll}
            & f(x_0) - f_\star \leq R^2 \\
            & \algo{\left< \nabla f(x_1), \nabla f(x_0) \right> = 0,} \\
            & \algo{\left< \nabla f(x_1), x_1 - x_0 \right> = 0.} \\
            & \forall i, j, ~ \class{f_i \geq f_j + \left< g_j, x_i - x_j \right> + \tfrac{1}{2L}\|g_i - g_j\|^2} \\
            & \class{\hspace{1.3cm} + \tfrac{\mu}{2(1 - \mu / L)}\|x_i - \tfrac{1}{L}g_i - x_j + \tfrac{1}{L}g_j\|^2.}
        \end{array}
        \right.
    \end{array}
    \right.
\end{equation*}

Using SDP lifting, we can formulate this problem as a semi-definite program of the form~\eqref{eq:pep_primal} using the variables
\begin{align*}
    F & = (f_\star, f_0, f_1)^{\top} \\
    G & = (x_\star, x_0, g_0, x_1, g_1)^{\top}(x_\star, x_0, g_0, x_1, g_1).
\end{align*}
We therefore set the parameters of~\eqref{eq:generic_proof} to the following values:

\scriptsize
\begin{align*}
    v_P & = (-1, 0, 1)^{\top}, \qquad M_P = 0_5, \\
    v_I & = (-1, 1, 0)^{\top}, \qquad M_I = 0_5,
\end{align*}
\begin{align*}    
    \class{v_{\mathcal{F}}^{(\star, 0)}} & = \begin{pmatrix}
        -1 \\ 1 \\ 0
    \end{pmatrix}, \quad
    \class{M_{\mathcal{F}}^{(\star, 0)}} = \frac{1}{2(1-\kappa)} \begin{pmatrix}
        \mu & -\mu & 1 & 0 & 0 \\
        -\mu & \mu & -1 & 0 & 0 \\
        1 & -1 & \frac{1}{L} & 0 & 0 \\
        0 & 0 & 0 & 0 & 0 \\
        0 & 0 & 0 & 0 & 0
    \end{pmatrix}
, \\
    \class{v_{\mathcal{F}}^{(\star, 1)}} & = \begin{pmatrix}
        -1 \\ 0 \\ 1
    \end{pmatrix}, \quad
    \class{M_{\mathcal{F}}^{(\star, 1)}} = \frac{1}{2(1-\kappa)} \begin{pmatrix}
        \mu & 0 & 0 & -\mu & 1 \\
        0 & 0 & 0 & 0 & 0 \\
        0 & 0 & 0 & 0 & 0 \\
        -\mu & 0 & 0 & \mu & -1 \\
        1 & 0 & 0 & -1 & \frac{1}{L}
    \end{pmatrix}, \\
    \class{v_{\mathcal{F}}^{(0, \star)}} & = \begin{pmatrix}
        1 \\ -1 \\ 0
    \end{pmatrix}, \quad
    \class{M_{\mathcal{F}}^{(0, \star)}} = \frac{1}{2(1-\kappa)} \begin{pmatrix}
        \mu & -\mu & \kappa & 0 & 0 \\
        -\mu & \mu & -\kappa & 0 & 0 \\
        \kappa & -\kappa & \frac{1}{L} & 0 & 0 \\
        0 & 0 & 0 & 0 & 0 \\
        0 & 0 & 0 & 0 & 0
    \end{pmatrix}, \\
    \class{v_{\mathcal{F}}^{(0, 1)}} & = \begin{pmatrix}
        0 \\ -1 \\ 1
    \end{pmatrix}, \quad
    \class{M_{\mathcal{F}}^{(0, 1)}} = \frac{1}{2(1-\kappa)} \begin{pmatrix}
        0 & 0 & 0 & 0 & 0 \\
        0 & \mu & -\kappa & -\mu & 1 \\
        0 & -\kappa & \frac{1}{L} & \kappa & -\frac{1}{L} \\
        0 & -\mu & \kappa & \mu & -1 \\
        0 & 1 & -\frac{1}{L} & -1 & \frac{1}{L}
    \end{pmatrix}, \\
    \class{v_{\mathcal{F}}^{(1, \star)}} & = \begin{pmatrix}
        1 \\ 0 \\ -1
    \end{pmatrix}, \quad
    \class{M_{\mathcal{F}}^{(1, \star)}} = \frac{1}{2(1-\kappa)} \begin{pmatrix}
        \mu & 0 & 0 & -\mu & \kappa \\
        0 & 0 & 0 & 0 & 0 \\
        0 & 0 & 0 & 0 & 0 \\
        -\mu & 0 & 0 & \mu & -\kappa \\
        \kappa & 0 & 0 & -\kappa & \frac{1}{L}
    \end{pmatrix}, \\
    \class{v_{\mathcal{F}}^{(1, 0)}} & = \begin{pmatrix}
        0 \\ 1 \\ -1
    \end{pmatrix}, \quad
    \class{M_{\mathcal{F}}^{(1, 0)}} = \frac{1}{2(1-\kappa)} \begin{pmatrix}
        0 & 0 & 0 & 0 & 0 \\
        0 & \mu & -1 & -\mu & \kappa \\
        0 & -1 & \frac{1}{L} & 1 & -\frac{1}{L} \\
        0 & -\mu & 1 & \mu & -\kappa \\
        0 & \kappa & -\frac{1}{L} & -\kappa & \frac{1}{L}
    \end{pmatrix},
\end{align*}

\begin{align*}
    \algo{v_{\mathcal{A}}^{(1)}} & = \begin{pmatrix}
        0 \\ 0 \\ 0
    \end{pmatrix}, \quad
    \algo{M_{\mathcal{A}}^{(1)}} = \frac{1}{2} \begin{pmatrix}
        0 & 0 & 0 & 0 & 0 \\
        0 & 0 & 0 & 0 & 0 \\
        0 & 0 & 0 & 0 & 1 \\
        0 & 0 & 0 & 0 & 0 \\
        0 & 0 & 1 & 0 & 0
    \end{pmatrix}, \\
    \algo{v_{\mathcal{A}}^{(2)}} & = \begin{pmatrix}
        0 \\ 0 \\ 0
    \end{pmatrix}, \quad
    \algo{M_{\mathcal{A}}^{(2)}} = \frac{1}{2} \begin{pmatrix}
        0 & 0 & 0 & 0 & 0 \\
        0 & 0 & 0 & 0 & -1 \\
        0 & 0 & 0 & 0 & 0 \\
        0 & 0 & 0 & 0 & 1 \\
        0 & -1 & 0 & 1 & 0
    \end{pmatrix}.
\end{align*}
\normalsize

Solving this SDP, we find the rate $\tau=\left(\frac{L-\mu}{L+\mu}\right)^2$.

Moreover, the corresponding dual values are

\setlength{\fboxrule}{1.5pt}
\begin{align*}
    & \fcolorbox{cyan!40}{white}{$\class{\lambda_\mathcal{F}^{\star, 0}} = \frac{2\mu(L-\mu)}{(L+\mu)^2},$} && \quad
    \fcolorbox{cyan!40}{white}{$\class{\lambda_\mathcal{F}^{\star, 1}} = \frac{2\mu}{L+\mu},$} \\
    & \class{\lambda_\mathcal{F}^{0, \star}} = 0, && \quad
    \fcolorbox{cyan!40}{white}{$\class{\lambda_\mathcal{F}^{0, 1}} = \frac{L-\mu}{L+\mu},$} \\
    & \class{\lambda_\mathcal{F}^{1, \star}} = 0, && \quad
    \class{\lambda_\mathcal{F}^{1, 0}} = 0, \\
    & \fcolorbox{purple!40}{white}{$\algo{\lambda_\mathcal{A}^1} = \frac{2}{L+\mu},$} && \quad
    \fcolorbox{purple!40}{white}{$\algo{\lambda_\mathcal{A}^2} = 1.$}
\end{align*}

Plugging those values in~\eqref{eq:generic_proof} builds a proof of convergence of~\eqref{eq:gdls} with the guarantee $f(x_1) - f_\star \leq \left(\frac{L-\mu}{L+\mu}\right)^2 (f_0 - f_\star)$.

~

\noindent
\resizebox{\linewidth}{!}{
\fbox{\parbox{\textwidth}{
    \vspace{-.4cm}
    \begin{align*}
        & f(x_1) - f_\star - \left(\frac{L-\mu}{L+\mu}\right)^2 (f(x_0) - f_\star) \\
        \leq &
        \fcolorbox{cyan!40}{white}{$
        \begin{array}{l}
            \frac{2\mu(L-\mu)}{(L+\mu)^2} \bigg( f(x_0) - f_\star + \left< \nabla f(x_0), x_\star - x_0 \right> + \tfrac{1}{2L}\|\nabla f(x_0)\|^2 + \tfrac{\mu}{2(1 - \mu / L)}\|x_\star - x_0 + \tfrac{1}{L}\nabla f(x_0)\|^2 \bigg) \\
            + \frac{2\mu}{L+\mu} \bigg( f(x_1) - f_\star + \left< \nabla f(x_1), x_\star - x_1 \right> + \tfrac{1}{2L}\|\nabla f(x_1)\|^2 + \tfrac{\mu}{2(1 - \mu / L)}\|x_\star - x_1 + \tfrac{1}{L}\nabla f(x_1)\|^2 \bigg) \\
            + \frac{L-\mu}{L+\mu} \bigg( f(x_1) - f(x_0) + \left< \nabla f(x_1), x_0 - x_1 \right> + \tfrac{1}{2L}\|\nabla f(x_0) - \nabla f(x_1)\|^2 + \tfrac{\mu}{2(1 - \mu / L)}\|x_0 - \tfrac{1}{L}\nabla f(x_0) - x_1 + \tfrac{1}{L}\nabla f(x_1)\|^2 \bigg)
        \end{array}
        $}
        \\
        &
        \fcolorbox{purple!40}{white}{$
        \begin{array}{l}
            + \frac{2}{L+\mu} \left< \nabla f(x_1), \nabla f(x_0) \right> \\
            + \left< \nabla f(x_1), x_1 - x_0 \right>
        \end{array}
        $}
        \\
        \leq & ~ 0.
    \end{align*}
    \vspace{-.8cm}
}}}

~

The first inequality holds independently on the chosen class. It simply results from terms rearrangement. By subtracting the LHS from the RHS, one would find a semi-definite positive quadratic form of the variables $x_0, x_1, \nabla f(x_0)$ and $\nabla f(x_1)$.
The second inequality precisely uses the (in)equalities that are specific to the chosen class and algorithm. Note that the two algorithm constraints can be replaced by the sole constraint $\left< \nabla f(x_1), x_1 - x_0 + \frac{2}{L+\mu} \nabla f(x_0) \right> = 0$, immediately showing that this guarantee also holds on the gradient descent method with fixed steps-size $\frac{2}{L+\mu}$.

\section{Lyapunov with PEPs}\label{sec:lyap}

We saw in \Cref{subsec:proofs_with_pep} that worst-case proofs essentially writes as~\eqref{eq:generic_proof}:

\vspace{-.6cm}

\begin{align}
    & \underbrace{\left<F, v_{P}\right> + \left<G, M_{P}\right>}_{\text{Performance metric}} - \tau \underbrace{\left[ \left<F, v_{I}\right> + \left<G, M_{I}\right> \right]}_{\text{Initialization}} \nonumber \\
    & \leq \sum_{j} \lambda^{(j)} \underbrace{\left[ \left<F, v^{(j)}\right> + \left<G, M^{(j)}\right> \right]}_{\text{Constraint}} \nonumber \\
    & \leq 0.
\end{align}

Namely, the right linear combination of the available constraints upper bounds the difference between the performance metric and $\tau$ times the initial value.
Sometimes those proofs can be relatively complicated and a simpler one can be desirable. In particular, this is the case when the algorithm under consideration is run for a few iterations.
Lyapunov analyses typically allows reducing the worst-case analyses of $T$ iterations to that of a single iteration, and therefore reducing the complexity of the proof.

For example, for~\eqref{eq:nag}, described as follow

\begin{equation}
    \begin{array}{c}
        \lambda_{t+1} = \frac{1}{2} + \sqrt{\frac{1}{4} + \lambda_t^2} \\
        y_t = x_t + \frac{\lambda_t - 1}{\lambda_{t+1}} (x_t - x_{t-1}), \\
        x_{t+1} = y_t - \frac{1}{L}\nabla f(y_t).
    \end{array}
    \tag{NAG}\label{eq:nag}
\end{equation}

on $\class{\F_{0, L}}$, we often use the sequence

\vspace{-.5cm}

\begin{equation}
    V_t = \lambda_{t}^2(f_t - f_\star) + \frac{L}{2}\|\lambda_{t}(x_t - x_\star) + (1-\lambda_{t})(x_{t-1} - x_\star)\|^2 \label{eq:naive_lyap}
\end{equation}
providing a worst-case convergence guarantee $f(x_T) - f_\star = \mathcal{O}(1/T^2)$.
In general, a direct way to find such a sequence is to consider

\vspace{-.5cm}

\begin{align}
    V_t = & \underbrace{\left[ \left<F, v_{I}\right> + \left<G, M_{I}\right> \right]}_{\text{Initialization}} \nonumber \\
    & + \underset{\text{before step} ~ t}{\underset{\text{values observed} }{\sum_{j ~ | ~ \text{only involves}}}} \lambda^{(j)} \underbrace{\left[ \left<F, v^{(j)}\right> + \left<G, M^{(j)}\right> \right]}_{\text{Constraint}}.
\end{align}

Applying this method on~\eqref{eq:nag} provides the sequence of complete potential functions

\vspace{-.5cm}

\begin{align*}
    V_t = & ~ \lambda_{t}^2(f_t - f_\star) + \frac{L}{2}\|\lambda_{t}(x_t - x_\star) + (1-\lambda_{t})(x_{t-1} - x_\star)\|^2 \\
    & ~ + \frac{1}{2L}\sum_{s=1}^{t-1} [ \lambda_{s+1}^2 \|\nabla f(x_{s+1})\|^2 + \lambda_{s+1} \|\nabla f(y_s)\| \\
    & \quad\quad\quad\quad\quad + \lambda_{s}^2\|\nabla f(y_s) - \nabla f(x_s)\|^2 ]
\end{align*}
that allows for free (using the same inequalities as for proving that \eqref{eq:naive_lyap} is decreasing) to also conclude that $\min_{t \leq T}\|\nabla f(x_t)\|^2 = \mathcal{O}(1 / T^3)$, as shown in~\cite[Theorem 5.2.d]{monteiro2013accelerated} and experimentally evidenced using PEPs in~\cite[Table 4]{taylor2017smooth}.

Note that the cumulatively summed up constraints involve both class constraints and algorithm constraints. Therefore, this technique can be applied directly on~\eqref{eq:gfom} while looking for an optimal algorithm, its rate, the corresponding proof and a sequence of potential functions at the same time.

\section{Conclusion}

\paragraph{Summary} not only is the \emph{performance estimation problem (PEP)} framework a powerful tool to automate the search of guarantees, but also it allows exhibiting general structure of proofs. Understanding this structure enables to generalize results onto larger class of functions or onto a class of methods, but also to find new optimization methods and study their convergence properties. Finally, it also enables to understand how to build a Lyapunov sequence of functions.

\paragraph{Open research directions} all this framework relies on two major assumptions: the class constraints are known and homogeneous in $\|x\|^2$ and $\|\nabla f(x)\|^2$ and $f$, and the method's update is a linear combination of previous iterates and observed oracle calls. Therefore, two interesting questions arise: can we automate the search of the interpolation conditions?
And, how can we generalize this framework to non homogeneous class of functions or to non linear methods such as adaptive step-size based methods? A few works already investigate this direction for some specific methods. In particular,~\cite{barre2020complexity} studies a variant of the Heavy-ball method~\cite{polyak_gradient_1963} using Polyak step-sizes, also discussed in~\cite[Chapter 4]{barre2021worst}. On the other hand,~\cite{gupta2023nonlinear} uses PEP techniques to provide worst-case guarantees on several variants of non-linear conjugate gradient methods.

\section*{Acknowledgments}

    The work of B. Goujaud and A. Dieuleveut is partially supported by ANR-19-CHIA-0002-01/chaire SCAI, and Hi!Paris.
    A.~Taylor acknowledges support from the European Research Council (grant SEQUOIA 724063).
    This work was partly funded by the French government under management of Agence Nationale de la Recherche as part of the ``Investissements d’avenir'' program,
    reference ANR-19-P3IA-0001 (PRAIRIE 3IA Institute).

\bibliographystyle{abbrv}
\bibliography{references}

\newpage

\end{document}